\pdfoutput=1
\documentclass[letterpaper, 10 pt, conference]{ieeeconf}  

\IEEEoverridecommandlockouts                              

\usepackage{graphicx} 
\usepackage{amsmath} 
\usepackage{amssymb}  
\usepackage{algorithm}                      
\usepackage{todonotes}
\usepackage[noend]{algpseudocode}
\usepackage{booktabs}
\newtheorem{remark}{Remark}
\title{\LARGE \bf
On Imitation Learning of Linear Control Policies:\\ Enforcing Stability and Robustness Constraints via LMI Conditions
}

\author{Aaron Havens and Bin Hu
\thanks{A. Havens, and B. Hu are with the Coordinated Science Laboratory (CSL) and the Department of Electrical and Computer Engineering,
        University of Illinois at Urbana-Champaign, Urbana
        {\tt\small \{ahavens2, binhu7\}@illinois.edu}}%
}

\begin{document}

\maketitle
\thispagestyle{empty}
\pagestyle{empty}

\begin{abstract}

When applying imitation learning techniques to fit a policy from expert demonstrations, one can take advantage of prior stability/robustness assumptions on the expert's policy and incorporate such control-theoretic prior knowledge explicitly into the learning process. In this paper, we formulate the imitation learning of linear policies as a constrained optimization problem, and present efficient methods which can be used to enforce stability and robustness constraints during the learning processes. Specifically, we show that one can guarantee the closed-loop stability and robustness by posing linear matrix inequality (LMI) constraints on the fitted policy. Then both the projected gradient descent method and the alternating direction method of multipliers (ADMM) method can be applied to solve the resulting constrained policy fitting problem. Finally, we provide numerical results to demonstrate the effectiveness of our methods in producing linear polices with  various stability and robustness guarantees.

\end{abstract}

\section{INTRODUCTION}

Recently, imitation learning (IL) for control design purposes
has received increasing  attention~\cite{argall2009survey,hussein2017imitation,osa2018algorithmic}. 
IL methods are particularly attractive for some control applications involving complex objective functions which are hard to specify beforehand. If we consider the control problem of autonomous driving, we must carefully trade off a few competing objectives such as performance, safety, and comfort level in driving. It can be difficult to come up with a precise metric for ``comfortable driving'', therefore, how to formulate the control design of self-driving as an optimization problem is unclear in the first place. It is known that human experts are rather adept at performing tasks such as driving and grasping. This provides a strong motivation for using IL methods which are designed to generate control policies that mimic human expert demonstration.

The basic idea of IL for control is to utilize expert demonstrations for policy fitting.
The simplest IL algorithm is the behavior cloning method which takes in demonstrations, such as a human driver's steering and throttle commands, and attempts to fit a state-action mapping in a supervised learning fashion. 
Such an approach typically requires a large amount of expert demonstrations.
Some more advanced IL algorithms such as DAgger~\cite{ross2011reduction} and GAIL \cite{ho2016generative} have been developed to handle issues such as covariate shift and sample efficiency. However, these methods do not generally provide stability and robustness guarantees when applied to control design. This  significantly restricts the practical deployment of these methods in real-world safety-critical applications. 



In contrast, many types of guarantees in terms of stability and robustness can be provided by modern control-theoretic methods.
Such guarantees have played a fundamental role in designing modern control systems~\cite{zhou1998essentials}, and can potentially be used as prior design knowledge in controller synthesis without necessarily constraining the control objective. It becomes natural to  ask how to incorporate these control-theoretic guarantees into the IL framework. 

In this paper, we study how to enforce stability and robustness guarantees for imitation learning of linear control policies.
We focus on the linear system case and view this as a benchmark for general IL. We present a constrained optimization formulation for IL in this setting, and show that various stability/robustness constraints in the form of linear matrix inequalities (LMIs) can be enforced during the policy fitting processes.
The resulting constrained IL problem has a special structure such that it can be efficiently solved using the projected gradient descent method or the alternating direction method of multipliers (ADMM). A comprehensive numerical study is also provided to demonstrate the effectiveness of our methods in producing linear policies with various stability/robustness guarantees.

%
%




%
%
%

\textbf{Related work:} Our work is mostly inspired by the recent result on linear policy fitting with Kalman constraints~\cite{kalman_constraint}.
Specifically, Malayandi et al. \cite{kalman_constraint} has proposed a method of learning linear policies from demonstrations which ensures that the policy is optimal with respect to some linear quadratic regulator (LQR) problem.
 However, the expert demonstrations need not be optimal in the LQR sense in the first place, and hence this approach may bias the policy fitting process for certain applications. In this paper, we address this bias issue by considering weaker prior assumptions such as stability or $\mathcal{H}_\infty$-robustness constraints. After the initial version of our paper was submitted to ACC 2021, we became aware of two concurrent papers~\cite{yin2020imitation, tu2021closing} which address IL with stability guarantees independently. 
Yin et al.~\cite{yin2020imitation} consider the design of neural network policies for linear time-invariant (LTI) systems and
propose an LMI-based IL method with stability guarantees in that setting. Tu et al.~\cite{tu2021closing}
propose the CMILe method which is capable of training nonlinear policies with the same safety guarantees as the experts. 
Our results complement these two papers by presenting a unified LMI-based treatment of stability and robustness constraints in IL of linear policies. There also exist some results on integrating the $\mathcal{H}_\infty$-robustness constraints into the reinforcement learning framework \cite{wu2013simultaneous, luo2014off,han2019h,zhang2019policy,zhang2020stability,zhang2021derivative, donti2020enforcing}. Our paper extends the use of
the $\mathcal{H}_\infty$-robustness constraint to the IL setting.


\section{PROBLEM FORMULATION}
\subsection{Notation}

The set of $n$-dimensional real vectors is denoted as
$\mathbb{R}^n$. The set of $m\times n$ real matrices is denoted as $\mathbb{R}^{m \times n}$.
The identity matrix is denoted as $I$. 
For brevity, repeated blocks in lower triangular parts of symmetric matrices are replaced 
by ``$*$''. When a matrix $P$ is positive semidefinite (definite), we will use the notation $P\succeq 0$ ($P\succ 0$). The spectral radius of $A$ is denoted as $\rho(A)$. The space $l_2$ is the space of square summable sequences with norm $||x||_2 = (\sum^{\infty}_i |x_i|^2)^{\frac{1}{2}}$. The Frobenius norm of a matrix $A \in \mathbb{R}^{n \times m}$ is denoted as $||A||_F$. 



\subsection{Background: Imitation Learning via Behavior Cloning}
One typical setting for IL is that we are given a control design problem whose objective function is hard to specify beforehand. 
The autonomous driving example mentioned in the introduction is exactly one such problem.
Then we can try to remedy this cost function design issue by directly fitting a policy on expert demonstrations. The hope is that the fitted policy can mimic the behaviors of experts and hence perform well on the given control problem.
 Specifically, suppose a sequence of state/action pairs $\{x_k, u_k\}_{k=0}^N$ has been demonstrated by the expert. Now we want to fit a function $u=K(x)$ from these observed demonstrations. Such a behavior cloning approach leads to the following finite-sum optimization problem:
 \begin{align}\label{eq:ERM}
     \underset{ K {}}{\text{minimize}} \quad \frac{1}{N}\sum^{N}_{k=0} l(K(x_k),u_k) + r(K)
 \end{align}
 where $l$ is some loss function measuring the empirical performance of the fitted policy on observed demonstrations, and $r$ is a regularization term introduced to prevent overfitting. In the above formulation, the policy is typically parameterized as a linear function or a neural network. 
 The resulting optimization problem is unconstrained and can be efficiently solved by stochastic gradient descent (SGD) or other first-order methods.
 Such a formulation can be problematic and inefficient for control applications 
where closed-loop stability, robustness, and safety become important 
 design concerns. For example, it may require a large number of demonstrations to ensure the solution of the unconstrained optimization problem \eqref{eq:ERM} to be a stabilizing policy, causing serious concerns on the efficiency and safety of the IL framework. Hence, it is more natural to adopt the following constrained optimization formulation:
  \begin{align}\label{eq:constrainedERM}
     \underset{ K\in\mathcal{K} {}}{\text{minimize}} \quad \frac{1}{N}\sum^{N}_{k=0} l(K(x_k),u_k) + r(K)
 \end{align}
where the policy fitting is confined to a feasible set $\mathcal{K}$ which is specified to carry the information of potential stability/robustness/safety constraints on $K$.
Such a constrained optimization formulation helps reduce sample complexity and improve system safety at the price of introducing more challenging computational tasks.  
In general, the set $\mathcal{K}$ is non-convex, and the constrained optimization problem \eqref{eq:constrainedERM} is more difficult to solve compared with its unconstrained counterpart \eqref{eq:ERM}. This motivates the study in this paper. Specifically, we will consider a simpler benchmark problem where the plant dynamics and the underlying expert policy are both assumed to be linear. We will discuss how to solve \eqref{eq:constrainedERM} with various stability and robustness constraints in this setting.

\subsection{Problem Setup: Constrained Linear Policy Fitting}
Now we confine the scope of our paper to linear policy fitting.  Suppose we want to fit a policy for a linear system with state $x_t$ and action $u_t$.
We assume that the ground truth expert policy is linear, i.e. \begin{align}\label{eq:expertK} 
u_t=K^* x_t+e_t
\end{align}
where $e_t$ is some zero mean noise. 
Notice that the appearance of $e_t$ is quite natural and intuitive since there will always be some noise in the human expert's demonstrations\footnote{In other words, humans will not make identical actions even given the same state observation.}. 
Based on these assumptions, it suffices to only consider linear state-feedback policy parameterized by a static gain matrix $K$. In general, it will be insufficient to use linear policy to model human experts' behaviors. It is our hope that our study on the simplified linear generative model \eqref{eq:expertK} can bring some insights for more realistic settings.

Suppose we have gathered the demonstrated state/action pairs $\{x_k, u_k\}_{k=0}^N$ from the generative model~\eqref{eq:expertK}. Here we use ``k" as the subscript to imply that the demonstrations may not be generated by a single trajectory of the underlying dynamical system.
 Our goal is to learn $K^*$ from the demonstrations. 
 We can adopt the loss function and the regularization term used in~\cite{kalman_constraint} to formulate the following optimization problem
\begin{align}\label{eq:imitate}
\underset{K \in \mathcal{K}{}}{\text{minimize}} \quad \frac{1}{N}\sum^{N}_{k=0} l(Kx_k,u_k)+r(K).
\end{align}
Typically, we can just set $l(Kx_k,u_k)=||  K x_k - u_k||_2^2$ and $r(K)=||K||_F^2$.
The feasible set $\mathcal{K}$ needs to be further specified. 
A naive option is to consider an unconstrained setting. However, prior knowledge of $K^*$ can be used to confine the policy search to a much smaller feasible set and improve the performance of IL with fewer demonstrations.
For example, the Kalman constraint can be enforced if $K^*$ is known to be the solution for some LQR problem~\cite{kalman_constraint}, and ADMM can be applied to solve the resulting constrained optimization problem efficiently \cite{kalman_constraint}. 
The Kalman constraint approach can significantly improve the sample efficiency of IL when the ground truth policy $K^*$ happens to be the solution for some LQR problem.
Of course, 
the assumption that $K^*$ satisfies the Kalman constraint can be too strong and introduce unnecessary biases into the policy search when $K^*$ is only sub-optimal in the LQR sense. For such scenarios, weaker constraints are preferred.
In this paper, we are mainly interested in the following two types of constraints.

\begin{enumerate}
\item Stability: Consider the following LTI system
\begin{align}\label{eq:lti_system}
x_{t+1} &= A x_t + B u_t+w_t
\end{align}
where $w_t$ is some stochastic process noise. Suppose $A \in \mathbb{R}^{n_x \times n_x}$, and $B \in \mathbb{R}^{n_u, n_u}$. 
A useful prior assumption on $K^*$ is that such that it should at least stabilize the closed-loop dynamics.
Therefore, to incorporate such prior knowledge, we can specify the feasibility set in \eqref{eq:imitate} as follows
\begin{align}\label{eq:feasi_stable}
    \mathcal{K}=\{K:\rho(A+BK)<1\}.
\end{align}
The stability constraint $\rho(A+BK)<1$ is much weaker than the Kalman constraint and will not introduce much bias into the policy fitting process.

\item $\mathcal{H}_\infty$-robustness:
Sometimes we will have prior knowledge on how robust $K^*$ is with respect to the model uncertainty. For example, we may know $K^*$ can stabilize the system \eqref{eq:lti_system} robustly in the presence of uncertainty in $A$ and $B$.
The following feedback interconnection model provides a general model for \eqref{eq:lti_system} subject to various types of uncertainties.
\begin{align}\label{eq:lti_expl}
\begin{split}
    x_{t+1} &= A x_t + B u_t+w_t+ B_1 v_t \\
    z_t &= C_1 x_t + D_{12} u_t \\
    v &=\Delta(z)
    \end{split}
\end{align}
where the uncertainty is modeled by a bounded operator $\Delta$ which is causal and maps any $\ell_2$ sequence $\{z_t\}$ to another $\ell_2$ sequence $\{v_t\}$. The above model is general since $\Delta$ can be set up properly to model uncertainty, nonlinearity, and time delays\footnote{For example, when the state/input matrices are not exactly known, the model \eqref{eq:lti_system} may be modified as
   $ x_{t+1}=(A+\Delta_A) x_t+(B+\Delta_B)u_t+w_t$,
which is a special case of \eqref{eq:lti_expl} with $B_1=I$ and $v_t=\Delta_A x_t+\Delta_B u_t$.}.
Based on the famous small gain theorem~\cite{zames1966input},
the system \eqref{eq:lti_expl}
is robustly stable for any $\Delta$ satisfying $\|\Delta\|_{\ell_2\rightarrow \ell_2}\le\frac{1}{\gamma}$ if we have $\rho(A+BK)<1$ and $||F(K)||_{\infty} < \gamma$ where $F(K)$ is an LTI system defined as
\begin{align*}
F(K) = 
\left( \begin{array}{c|c}
   A + B K & B_1 \\
   \midrule
   C_1 + D_{12}K & 0\\
\end{array}\right).
\end{align*}
Notice $||F(K)||_{\infty}$ denotes the 
$\mathcal{H}_\infty$-norm of $F(K)$.
Therefore, if we know that the ground truth expert policy $K^*$ achieves the $\gamma$-level robustness, i.e. $||F(K^*)||_{\infty}<\gamma$, we can enforce such a constraint during the policy fitting process  by specifying 
\begin{align}
\label{eq:feasi_Hinf}
\mathcal{K}=\{K:\rho(A+BK)<1\,\,,\,\,||F(K)||_{\infty}<\gamma\}.
\end{align}
\end{enumerate}
It is worth mentioning that 
there exist many other types of stability and robustness constraints which can be posed to confine the policy search. For readability and clarity, our paper focuses on the above two commonly-used constraints.

\begin{remark}
For the ease of exposition, we will mainly focus on enforcing stability/robustness constraints for behavior cloning. Similar ideas can be applied to enforce constraints for DAgger which was originally developed to address the
covariate shift issue. Intuitively, the covariate shift issue is less significant for linear policy fitting, and hence we will skip the detailed discussion on DAgger.
\end{remark}

\subsection{Enforcing Stability/Robustness Constraints via LMIs}

 The feasible set specified by \eqref{eq:feasi_stable} or \eqref{eq:feasi_Hinf}
 is not convex in $K$, causing trouble for the constrained policy optimization. 
Fortunately, these stability/robustness conditions can be convexified as LMIs if we change variables properly. Now we briefly review these LMI conditions below.

\begin{enumerate}
\item Convex conditions for stabilization: Notice that $K$ stabilizes the system \eqref{eq:lti_system} if and only if there exists 
a symmetric positive definite matrix $P \in \mathbb{R}^{n_x\times n_x}$ such that 
  $ P - (A+B K)^{\intercal} P (A+BK) \succ 0$. This condition is bilinear in $P$ and $K$, leading to a non-convex control design condition. However, we can apply the Schur complement lemma to obtain the following intermediate matrix inequality.
\begin{align*}
    \begin{bmatrix}
    P^{-1} & (A+BK)P^{-1}\\
    * & P^{-1}
    \end{bmatrix}
    \succ 0,\quad P^{-1} \succ 0
\end{align*}
Then we can change variables as $(Q,L)= (P^{-1}, K P^{-1})$ and obtain the following LMI condition:
\begin{align}\label{eq:lmi_1}
    \begin{bmatrix}
    Q & A Q + B L\\
    * & Q
    \end{bmatrix}
    \succ 0,\quad Q \succ 0
\end{align}
The above reparameterization is well-known~\cite{duan2013lmis}. Now we have a convex set of $(Q,L)$ which can be used to extract stabilizing policies efficiently.
Since there is a one-to-one correspondence between $(P,K)$ and $(Q, L)$, such a reparameterization does not introduce any conservatism. 
\item Convex conditions for $\mathcal{H}_\infty$-robustness: We can impose an LMI condition which is similar to ~\eqref{eq:lmi_1} to describe all $K$ lying in the set \eqref{eq:feasi_Hinf}. Specifically, we can set $K = L Q^{-1}$ to satisfy $\rho(A+BK)<1$ and $||F(K)||_{\infty} < \gamma$ if there exists $Q \in \mathbb{R}^{n_x\times n_x}$ and $L \in \mathbb{R}^{n_u\times n_x}$ such that the following LMI is satisfied.
\begin{align}\label{eq:lmi_robust}
    \begin{bmatrix}
        Q & A Q + B L & B_1 & 0\\
    * & Q & 0 &  Q C_1^{\intercal} + L^{\intercal} D_{12}^{\intercal}\\
    * & * & I & 0\\
    * & * & * & \gamma^2 I
    \end{bmatrix}
    \succ 0,\, Q \succ 0.
\end{align}
Again, the above condition is convex in $(Q,L)$ and can be useful for our constrained policy fitting problem. As $\gamma \rightarrow \infty$, the $\mathcal{H}_\infty$-robustness condition \eqref{eq:lmi_robust} reduces to \eqref{eq:lmi_1}. Therefore, the $\mathcal{H}_\infty$-robustness assumption is a stronger prior which may be useful when $K^*$ happens to be robust. However, it may introduce some bias into policy fitting if $K^*$ is not robust in the first place.  
\end{enumerate}

From the above discussions, it becomes obvious that we can change the variable in the constrained policy fitting problem \eqref{eq:imitate} to obtain a problem with a convex constraint. Notice such a parameterization method has been used in safe reinforcement learning \cite{friedrich2017robust}. Here we use it in the IL setting. To summarize, the constrained IL problem can be recast into the following general form: 
\begin{align}\label{eq:stable_imitate}
&& \underset{(Q, L)}{\text{minimize}} \quad &\frac{1}{N}\sum^{N}_{k=0} l(L Q^{-1} x_k , u_k)+r(LQ^{-1})\\
&& \text{subject to}\quad&\text{LMI}(Q,L)\notag
\end{align}
where $\text{LMI}(Q,L)$ is some LMI with decision variables $(Q,L)$. 
If the only prior assumption is that the expert's policy stabilizes the closed-loop dynamics, the stability constraint \eqref{eq:lmi_1} will be used. If we further believe that the expert has achieved some level of robustness, the $\mathcal{H}_\infty$-robustness constraint \eqref{eq:lmi_robust} can be applied with some tuned $\gamma$.
Although the objective function in \eqref{eq:stable_imitate} is non-convex, the constraint $\text{LMI}(Q,L)$ is convex and can be easily handled using projected gradient methods.

\begin{remark}
It is worth emphasizing that \eqref{eq:stable_imitate} provides a general formulation for constrained imitation learning of linear policies due to the existing large body of stability/robustness conditions in the form of $\text{LMI}(Q,L)$~\cite{boyd1994linear,duan2013lmis, caverly2019lmi}. 
For example, one can consider time-varying polytopic uncertainty and enforce a robust stability constraint in the form of $\text{LMI}(Q,L)$. In addition, one can even consider neural network policies, and the local stability constraint used in \cite{yin2020imitation} is exactly in the form of $\text{LMI}(Q,L)$. For the ease of exposition, our paper mainly focuses on the stability and $\mathcal{H}_\infty$-robustness constraints.
\end{remark}

\section{Main Algorithms}


In this section, we present two optimization methods for solving the reformulated constrained IL problem \eqref{eq:stable_imitate}.

\subsection{Projected Gradient Descent}




Although the objective function in~\eqref{eq:stable_imitate} is non-convex, the projection onto the feasible set of the decision variables $(Q,L)$ is easy to compute. Therefore, simple optimization methods such as the projected gradient method or projected SGD can be readily applied.
The projection step can be performed by solving the following convex subproblem.
\begin{align}\label{eq:projection}
&&\Pi_{\mathcal{K}}(\hat Q, \hat L) = \underset{(Q,L)}{\text{argmin}}\quad &||(\hat Q,\hat L) - (Q,L)||_F^2\\
&&\text{subject to}\quad  &\text{LMI}(Q,L)\notag
\end{align}

Now we are ready to present our main algorithm which can be used to efficiently solve the constrained IL problem~\eqref{eq:stable_imitate}. The method that we use is the projected gradient descent (PGD) method, which is summarized as follows.
\begin{algorithm}
\begin{algorithmic}
\caption{Projected Gradient Descent}
    \label{alg:PGD}
\Procedure{Projected Gradient Descent}{A,B}
\State $(Q^{(0)}, L^{(0)}) \gets \text{Find a feasible solution of LMI}(Q,L)$
\For{$n=1:T$ iterations}
    \State $(Q^{(n)}, L^{(n)}) \gets \Pi_{\mathcal{K}}(SGD(Q^{(n-1)},L^{(n-1)}))$
\EndFor
\State \textbf{return} $K = L^{(N)} (Q^{(N)})^{-1}$
\EndProcedure
\end{algorithmic}
\end{algorithm}

The matrices $(Q, L)$ are initialized by first finding a feasible solution of $\text{LMI}(Q,L)$. The gradient update performed on the objective in~\eqref{eq:stable_imitate} is performed with a first-order batch method such as SGD~\cite{ruder2016overview}. If the number of samples is small (i.e. less than $\approx 50$), we compute the gradient with a single batch which reduces to regular gradient descent. The gradient can be calculated using any auto-differentiation method. In the next section, we will choose to use PyTorch due to its efficient implementation of batch operations in conjunction with auto-differentiation. The convex projection problem~\eqref{eq:projection} and the feasibility of all LMIs can be  solved using the CVXPY modeling framework~\cite{diamond2016cvxpy} and the MOSEK solver~\cite{andersen2000mosek}.




\subsection{Alternating Direction Method of Multipliers (ADMM)}
By introducing the equality constraint $L=KQ$, we can obtain another useful equivalent form for \eqref{eq:imitate} as follows
\begin{align}\label{eq:stable_imitate_form2}
&& \underset{(K, Q, L)}{\text{minimize}} \quad &\frac{1}{N}\sum^{N}_{k=0} l(Kx_k , u_k)+r(K)\\
&& \text{subject to}\quad& L=KQ,\,\,\mbox{and}\,\,\text{LMI}(Q,L)\notag
\end{align}
The decision variables for the above optimization problem are $(K,Q,L)$. 
Notice that \eqref{eq:stable_imitate_form2} has a convex objective in $K$, a bilinear equality constraint in $(K, Q, L)$, and a convex inequality constraint in $(Q, L)$. This special structure allows us to apply ADMM just as in~\cite{kalman_constraint}.
Denote $C(K)=\frac{1}{N}\sum^{N}_{k=0} l(Kx_k , u_k)+r(K)$.
Then the augmented Lagrangian is given as
$\mathcal{L}_{\rho}(K,Q,L,Y) = C(K)+ \text{Tr}(Y^{\intercal}(K Q - L)) + \frac{\rho}{2} || K Q - L ||_F^2$,
where $Y$ is the dual variable corresponding to the constraint $KQ = L$ and $\rho$ is a fixed scalar penalty parameter which must be sufficiently large to yield convergence~\cite{gao2020admm} (we typically set $\rho=1$). Given initial parameters $(Q^{(0)},L^{(0)}, Y^{(0)})$, the updates of ADMM follow a \textit{three-step} optimization procedure, where each step is a convex problem:\\
 1. In the $K$-step, we update the variable $K$ as
\[            
K^{(n+1)} = \underset{K}{\text{argmin}}\quad \mathcal{L}_{\rho}(K, Q^{(n)}, L^{(n)}, Y^{(n)}).
\]
 2. In the $(Q, L)$-step, we update $(Q,L)$ by solving the following semidefinite program:
        \begin{align*}
            &&(Q^{(n+1)}, L^{(n+1)}) =\underset{(Q,L)}{\text{argmin}}\quad &\mathcal{L}_{\rho}(K^{(n+1)}, Q, L, Y^{(n)})\\
            && \text{subject to}\quad& \text{LMI}(Q,L)\notag
        \end{align*}
 3. In the $Y$-step, we update the variable $Y$ as
    \begin{align*}
        Y^{(n+1)} = Y^{(n)} + \rho (K^{(n+1)} Q^{(n+1)} - L^{(n+1)})
    \end{align*}

We can run the algorithm for a fixed number of iterations. Since the subproblems in ADMM are all convex, only CVXPY will be needed in the implementations. Note that ADMM is not necessarily superior to PGD in all settings, as we will see in the numerical case studies. 

\begin{remark}
Based on the results in \cite{yin2020imitation}, one can enforce local stability constraints for imitation learning of neural network policies using a constrained optimization formulation similar to \eqref{eq:stable_imitate_form2}. The only difference is that for neural network policies, the cost function becomes non-convex in $K$ and hence the $K$-step in ADMM can only be solved approximately using gradient-based methods. It is possible that robustness constraints on neural network policies can be enforced in a similar way.
\end{remark}


\section{NUMERICAL RESULTS}
In this section, we perform numerical studies to validate the effectiveness of the proposed methods. We mainly consider three different ground truth expert policies: i) aggressive stabilizing controllers;
    ii) LQR optimal controllers;
    and iii) $\mathcal{H}_{\infty}$ robust optimal controllers.
We will compare the proposed methods (PGD and ADMM) against the standard unconstrained policy fitting (PF) formulation~\eqref{eq:imitate} and the Kalman constraint method in~\cite{kalman_constraint} in terms of a defined cost appropriate for each setting. We generate the system trajectories by simulating the following model: 
\begin{align*}
    x_{t+1} &= A x_t + B u_t + w_t,\quad w_t \sim \mathcal{N}(0,W)\\
    u_t &= K^* x_t + e_t,\quad e_t \sim \mathcal{N}(0,\Sigma),\quad x_{0 i}\sim \mathcal{U}(-1,1),
\end{align*}
where $W$ and $\Sigma$ are the covariance matrices for $w_t$ and $e_t$, respectively. Here $\mathcal{N}$ denotes the normal distribution, and $\mathcal{U}$ denotes the uniform distribution.
The plant is generated by sampling state/input matrices $(A,B)$ from a uniform distribution. For illustrative purposes, the state dimension and control dimensions are fixed to be $n_x=4$ and $n_u=2$. We set $A_{ij},\, B_{ij} \sim \mathcal{\mathcal{U}}(-1, 1)$, and $\Sigma = W = (0.25) I$.
We keep on sampling until getting a stabilizable pair $(A,B)$.
\subsection{Aggressive Stabilizing Demonstrations}\label{aggressive}
\begin{figure}
    \centering
    \includegraphics[width=1.0\linewidth]{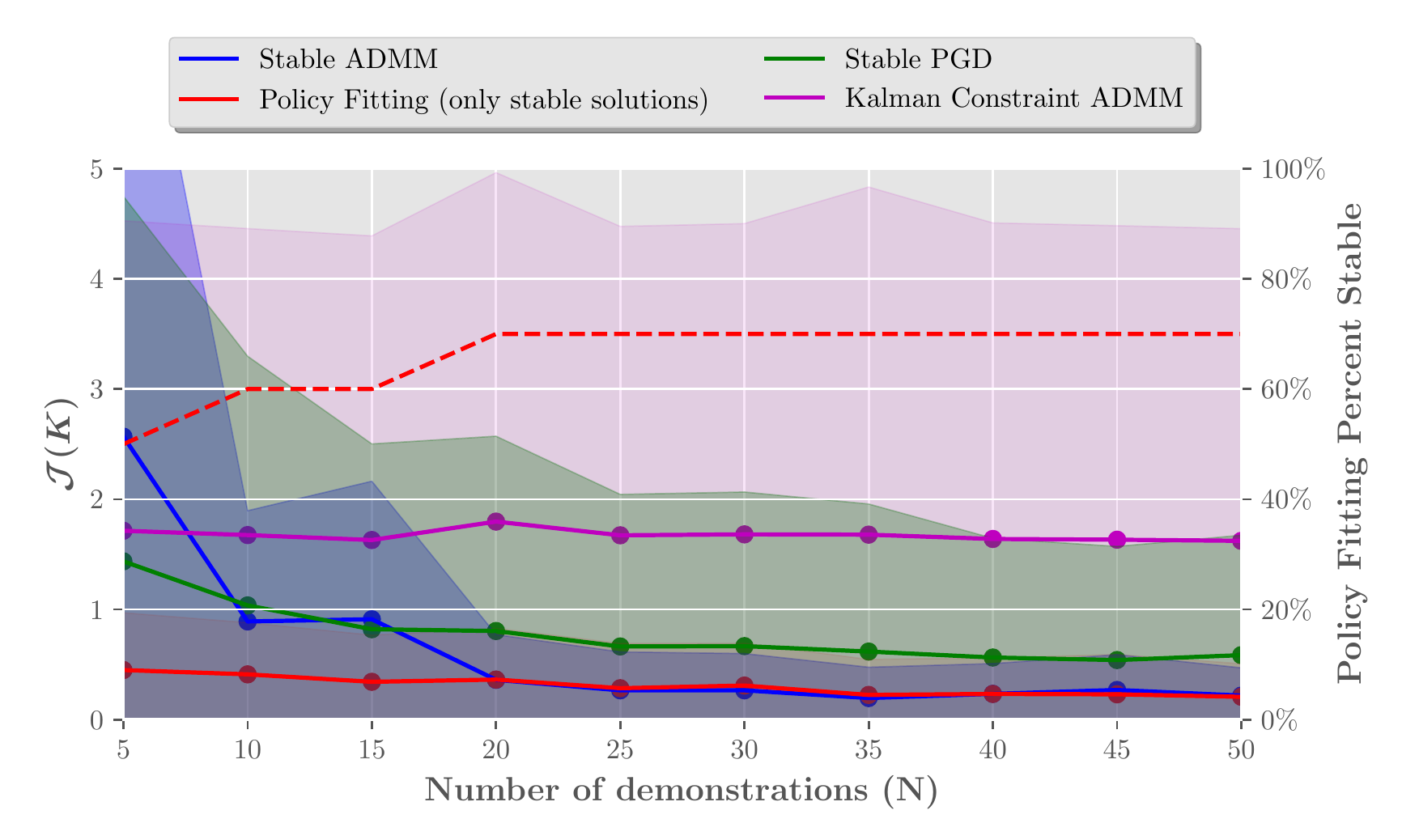}
    \caption{\textbf{Aggressive Stabilizing Demonstrations}: We compare the prediction performance of standard PF, stable PGD, stable ADMM and Kalman-constrain ADMM  with respect to the true demonstrator policy as the number of demonstrations $N$ increases. The standard PF problem frequently produces unstable controllers and where the constrained methods produce stable policies at every sample and data setting $N$. We sample $10$ random stabilizable systems with aggressive stabilizing controller$\{A_i, B_i, K_i\}_{i=1}^{10}$ and sample demonstrations in increments $N \in \{5,10\ldots,50\}$, showing the mean and $1$ standard deviation for each increment. It is important to note that for unconstrained PF we only display the mean of stable solutions and the red dashed line indicates the percentage of stable solutions.}
    \label{fig:stable_exp}
\end{figure}

To evaluate the application of our approach for generally stabilizing policies that are not necessarily robust, we consider expert policies which are aggressive and not very robust (e.g. $A$ has eigenvalues which are close to $1$). Specifically, we generate $K^*$ by sampling uniform random initializations of $(Q,L)$ and then projecting via the LMI \eqref{eq:lmi_1}.
Note that it is quite unlikely for the resultant matrix $K$ to be an optimal policy in either an LQR or $\mathcal{H}_{\infty}$ sense since any random initialization of $(Q, L)$ outside the stability feasible set will be projected to some interior point near the boundary of the feasible set. For each $K$, we measure how much the resulting state sequence $\{x_i\}_{i=1}^{N_t}$ under the closed-loop dynamics differs from the true closed-loop under $K^*$. Given the same initial condition, the cost is defined as 
$\mathcal{J}(K) = \frac{1}{N_{t}}\sum^{N_{t}}_{i=1} ||x_i^*-x_i||_2^2$. In our simulations, the cost is averaged over $100$ random test trajectories of length $N_t=100$. This cost is chosen over the regular PF objective in~\eqref{eq:imitate} since it reflects the actual system behavior over time. To represent the performance of standard PF more fairly, we only average over \textit{stable} solutions, while showing the percentage of stable solutions denoted by the red-dashed line in Figure~\ref{fig:stable_exp}. We show in Figure~\ref{fig:stable_exp} that the standard PF problem~\eqref{eq:imitate} frequently produces unstable policies where the stability-constrained methods of PGD and ADMM maintain a low cost at all times. On the other hand, the Kalman constraint is biased towards LQR solutions with sufficient robustness margins, and hence a large gap persists between the stability constrained methods in this setting. We also note that there is a gap between PGD and ADMM, which may be due to the hyperparameters chosen in our simulations.

\subsection{$LQR$ Optimal Demonstrations}
Next, we consider expert policies which are optimal with respect to the standard infinite horizon quadratic cost.
For evaluation, a finite, yet large, horizon of ($N_t = 1000$) is used to approximate the infinite horizon cost. We define the evaluation cost as 
    $\mathcal{J}(K) = J^{(N_t)}(K) - J^{(N_t)}(K^*)$.
For illustrative purposes, we set $Q = R = I$. Since the optimal policy for an LQR problem has guaranteed stability margins, the expert policy $K^*$ is well within the feasible set and standard PF does not fail as often as in the previous aggressive case. While all constrained methods remain stable, it is hard to beat the inductive bias given by the Kalman constraint in the LQR setting. With enough data, our stable ADMM method approaches the performance of the Kalman constraint method. In this case, stable PGD performs better than ADMM in the lower data regime, while there is a small gap as $N$ grows, possibly due to the hyperparameter choices.
\begin{figure}
    \centering
    \includegraphics[width=1.0\linewidth]{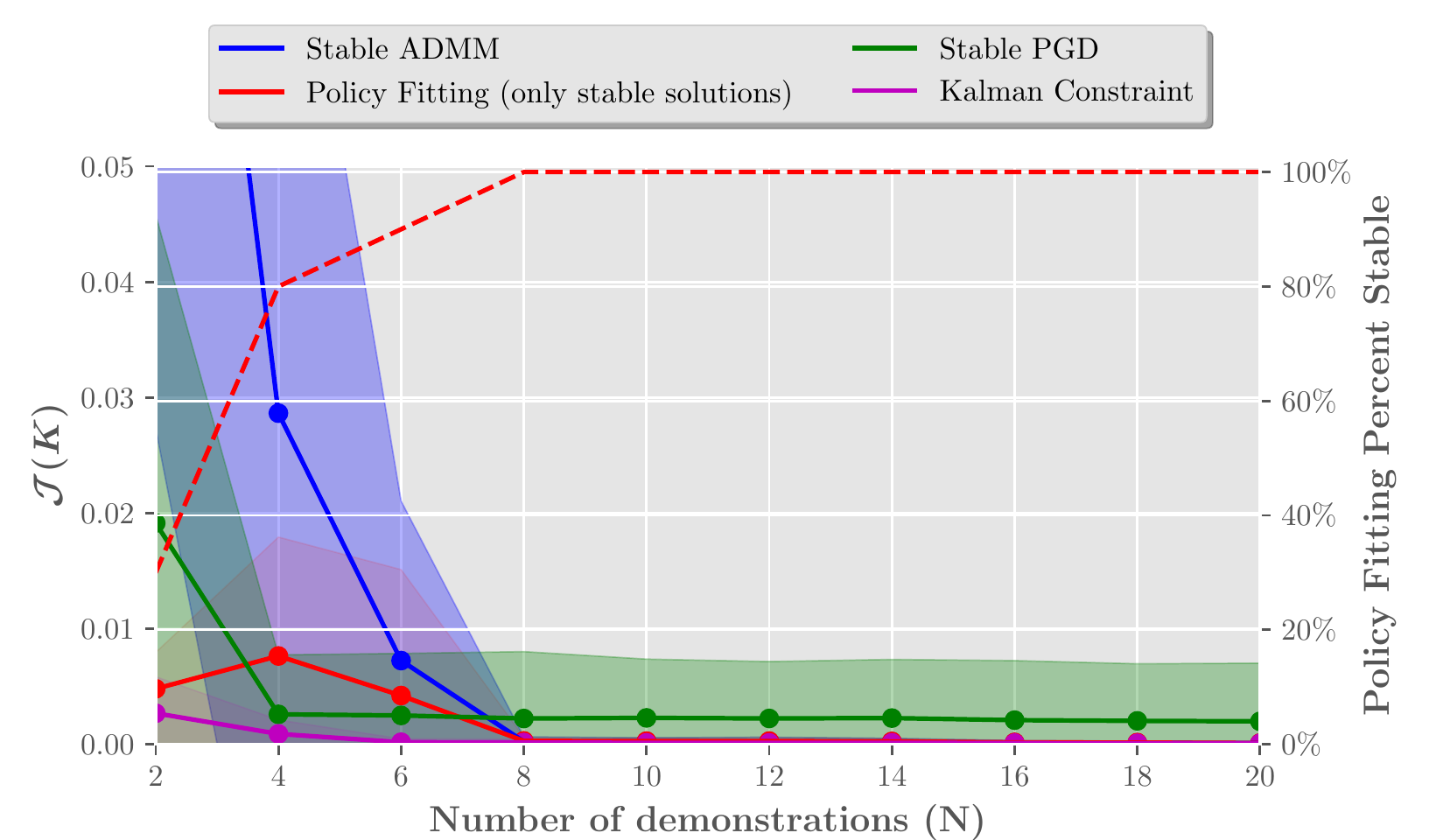}
    \caption{\textbf{LQR Demonstrations}: We sample $\{A_i, B_i, K_i\}_{i=1}^{10}$ and demonstrations in increments $N \in \{2,4\ldots,20\}$. We compare all constrained methods against standard PF with respect to the optimal LQR cost. Again, only stable solutions are averaged for standard PF and percentage of stable solutions is denoted by the red dashed line.}
    \label{fig:lqr_exp}
\end{figure}
\subsection{$\mathcal{H}_{\infty}$ Optimal Demonstrations}
Lastly, we explore a useful application of the $\mathcal{H}_\infty$-robustness constraint~\eqref{eq:lmi_robust} which ensures $||F( K)||_{\infty} < \gamma$ throughout the policy learning processes. For this experiment we choose $\gamma = \gamma^* + 0.1$ where $\gamma^*$ is the optimal $H_{\infty}$ norm solution. Suppose we have prior knowledge of the robustness present in the system and choose $\gamma$ as an upper-bound estimate of the true $\mathcal{H}_{\infty}$ norm. In this case the demonstrator is a noisy $\mathcal{H}_{\infty}$ optimal controller for the sampled stabilizable system where 
    $ A_{ij},\, B_{1ij}, B_{2ij} \sim \mathcal{\mathcal{U}}(-1, 1)$, $C_1 = I$, and $ D_{12 ij} \sim \mathcal{U}(-0.1,0.1)$.
Here we have $B_1 \in \mathbb{R}^{n_x \times n_d},D_{12} \in \mathbb{R}^{n_z \times n_d}$, $C_1 \in \mathbb{R}^{n_z,n_x}$ and $(n_x, n_u, n_d, n_z) = (4,2,1,4)$. The cost is given as
$\mathcal{J}(K) = ||F(K)||_{\infty} - ||F(K^*)||_{\infty}$.
To study the effect of this additional $\mathcal{H}_\infty$-robustness constraint, we compare against the first proposed stability constraint (using ADMM), standard PF and the Kalman constraint method. We see that standard PF still produces some unstable solutions, but all the constrained IL methods eventually perform well. Stable ADMM and the Kalman constraint method produce stabilizing controllers through each trial, but with a much larger $\gamma$ initially. The Kalman constraint has inherent built-in margins which explains its performance over stable-ADMM when there are few demonstrations available. The robust ADMM and PGD methods produce an $\mathcal{H}_{\infty}$ norm much lower than other methods and stay within the prescribed $\gamma$ threshold for all trials.

\begin{figure}
    \centering
    \includegraphics[width=1.\linewidth]{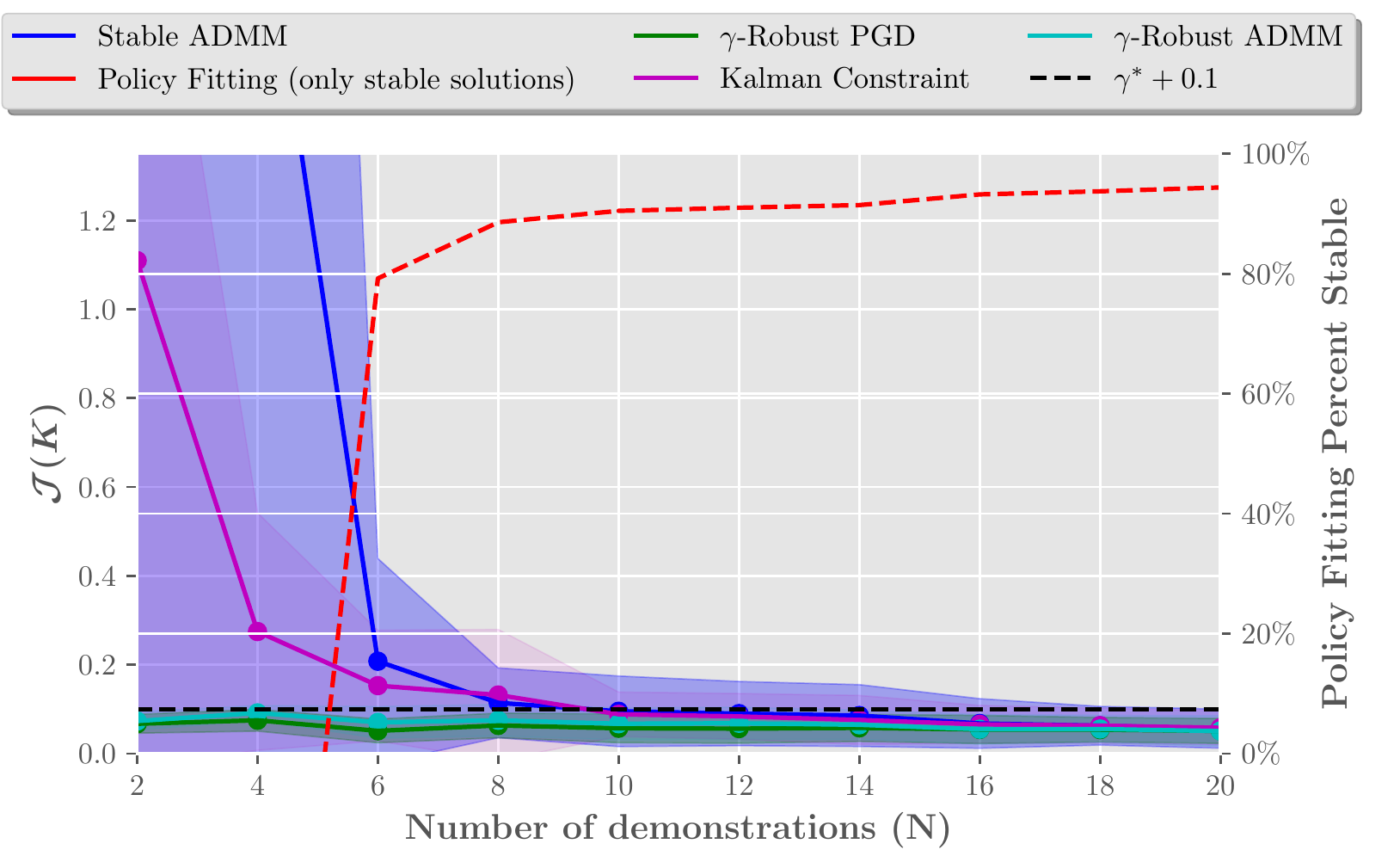}
    \caption{$\mathcal{H}_{\infty}$\textbf{-Optimal Demonstrations}: On a sample of $10$ systems, the $\gamma$-robust constrained methods against the previous stable ADMM, Kalman constraint and standard PF. We prescribe $\gamma = \gamma^* + 0.1$ for the robust projection as prior knowledge that the demonstrator is robust. We see that the robust constraint benefits even over the stable ADMM method and other baselines, while staying within the prescribed threshold.}
    \label{fig:plot_robust}
\end{figure}

\section*{ACKNOWLEDGMENT}
This work is generously supported by the NSF award CAREER-2048168.





\bibliography{bibtex}
\bibliographystyle{IEEEtran}

\end{document}